\documentclass[12pt]{article}
\usepackage{amsmath, amsthm, amscd, amsfonts, amssymb, graphicx, color}
\usepackage{maplestd2e}
\usepackage{tikz}
\usepackage{url}

\newtheorem{lemma}{Lemma}[section]
\newtheorem{remark}[lemma]{Remark}
\newtheorem{theorem}[lemma]{Theorem}

\newtheorem{example}[lemma]{Example}

\newcommand{\be}{\begin{equation}}
	\newcommand{\ee}{\end{equation}}
\begin{document}
	\setlength{\unitlength}{1mm}
	\baselineskip 8mm

\thispagestyle{empty}
\title{
On Gosper-Karaji algebraic Identities}
\author
{Hossein Teimoori Faal
	\vspace{.25in}\\
	Department of Mathematics and Computer Science,\\ 
	Allameh Tabataba'i University, Tehran, Iran
}

	\maketitle
	
\begin{abstract}

In this paper, we first quickly review the basics of an algebro-geometric method of Karaji's L-summing technique 
in today's modern language of algebra. 
Then, we also review the theory of Gosper's algorithm as 
a decision procedure for obtaining the indefinite sums 
involving \emph{hypergeometric} terms. 
Then, we show that how one can use Gosper's algorithm equipped with the L-summing method to obtain a class of combinatorial identities associated with a given algebraic identity. 
	
\end{abstract}

\section{Introduction}

About $1000$ years ago, an interesting algebro-geometric method was
used by M. Karaji \cite{WF}; an Iranian scholar, to give the first known
proof of the following well-known algebraic identity:
$$1^3+2^3+\cdots+n^3=(1+2+\cdots+n)^2.$$
Here we give the abstract idea behind his interesting approach, using today's
modern algebraic language. 
Consider a square array of numbers
$A_{n}=(a_{ij})_{1\leq i,j \leq n}$, as in table \ref{Table1}. 
Our main goal here is to
find the sum of all entries of the array 
$
A_{n}
$
in two different ways. We denote this sum 
by 
$
S_n=\sum_{i,j=1,\ldots,n} a_{ij}
$.
First, summing by rows gives us
$$
S_n=\sum_{k=1}^{n}\left(\sum_{i=1}^{n}a_{ki}\right).
$$

Then, summing by L-shape pieces (squared pieces in Table \ref{Karaji1}) results in
the following formula for $S_n$
$$
S_n=\sum_{k=1}^{n}\left(\sum_{i=1}^{k}a_{ik}+\sum_{j=1}^{k}a_{kj}-a_{kk}\right).
$$
Therefore, we get the following general identity:

\begin{eqnarray}\label{KarajiEq1}
	\sum_{k=1}^{n}\left(\sum_{i=1}^{n}a_{ki}\right)=\sum_{k=1}^{n}\left(\sum_{i=1}^{k}a_{ik}+\sum_{j=1}^{k}a_{kj}-a_{kk}\right).
\end{eqnarray}

\begin{center}
	\begin{table}\label{Table1}
		\centering
		\begin{tabular}{|c||ccccccc|} \hline
			$\downarrow i\setminus j\rightarrow$ &1& 2 & $\cdots$&$k$&$\cdots$&$n-1$&$n$\\
			\hline\hline
			$1$ &$a_{11}$&$a_{12}$&$\cdots$&$\framebox{$a_{1k}$}$&$\cdots$&$a_{1,n-1}$&$a_{1n}$\\
			$2$ &$a_{21}$&$a_{21}$&$\cdots$&$\framebox{$a_{2k}$}$&$\cdots$&$a_{2,n-1}$&$a_{2n}$\\
			$\vdots$ &$\vdots$&$\vdots$&$\ddots$&$\vdots$&$\vdots$&$\vdots$&$\vdots$\\
			$k$&$\framebox{$a_{k1}$}$&$\framebox{$a_{k2}$}$&$\cdots$&$\framebox{$a_{kk}$}$&$\cdots$&$a_{k,n-1}$&$a_{kn}$\\
			$\vdots$ &$\vdots$&$\vdots$&$\vdots$&$\vdots$&$\ddots$&$\vdots$&$\vdots$\\
			$n-1$ &$a_{n-1,1}$&$a_{n-1,2}$&$\cdots$&$a_{n-1,k}$&$\cdots$&$a_{n-1,n-1}$&$a_{n-1,n}$\\
			$n$ &$a_{n1}$&$a_{n2}$&$\cdots$&$a_{nk}$&$\cdots$&$a_{n-1,n}$&$a_{nn}$\\
			\hline
		\end{tabular}
		\caption{Karaji's L-summing method\label{Karaji1}}
	\end{table}
\end{center}

From here on, we will call the identity (\ref{KarajiEq1}) the \emph{L-summing equation}. As an
immediate consequence of the equation (\ref{KarajiEq1}), we have the following
identity

\begin{eqnarray}\label{KarajiSpecial1}
	\sum_{k=0}^{n}\left(\sum_{i=0}^{k}f(i)+k
	f(k)\right)=(n+1)\sum_{i=0}^{n}f(k),
\end{eqnarray}
where it is just a simple substitution of the entry $a_{ij}$ with a discrete function of one varaiable $f(i)$.
Next, we give several identities which can be derived
directly from formulas (\ref{KarajiEq1}) and (\ref{KarajiSpecial1}).
\\

\begin{example}[An identity for zeta function]
	
If we define $\zeta_{n}(s)=\sum_{k=1}^{n}\frac{1}{k^s}\hspace{2mm}(s
\in  \Re(s))$, then we get the zeta function
$\zeta(s)=\sum_{k=1}^{\infty}\frac{1}{k^s}=lim_{n\rightarrow\infty}\zeta_{n}(s).$
Now, let $a_{ij}=\frac{1}{i^s}\times\frac{1}{j^s}$ in formula (\ref{KarajiEq1}). 
Then, after applying the limit operator, we get
$$\sum_{k=1}^{\infty}\frac{\zeta_{k}(s)}{k^s}=\frac{\zeta^{2}(s)+\zeta(2s)}{2}.
$$
	
\end{example}

\begin{example}[An identity for derangements $D_{n}$]
	
Recall that the derangment $D_{n}$ counts the number of permutations of the set 
$[n]=\{1,2,\ldots,n\}$ 
without \emph{fixed-points}. 
Let $\Gamma(a,z)=\int_{z}^{\infty}e^{-t}t^{a-1}dt$. 
Then we have
$D_{n}=e^{-1}\Gamma(n+1,-1)$(see \cite{H1}). 
Now, by choosing
$f(i)=\Gamma(n+1,-1)=e\times D_{n}$ in formula (\ref{KarajiSpecial1}), we obtain the
following combinatorial identity for derangements
	
$$
\sum_{k=1}^{n}\left(\sum_{i=1}^{k}D(i)+k
D(k)-D(k)\right)=n\sum_{i=1}^{n}D(k).
$$ 
	
\end{example}
For more interesting examples see the reference \cite{H2}.
\\
Using similar arguments, based on the three dimensional geometric
idea and the principle of \emph{inclusion-exclusion}, we can obtain a 
generalization of Karaji's L-summing method for three
dimensional array $A_{n}=(a_{ijk})_{1\leq i,j,k \leq n}$.

\begin{eqnarray}\label{ThreeDimKaraji}
	&&\sum_{k=1}^{n}\sum_{j=1}^{n}\sum_{i=1}^{n}a_{ijk}\nonumber\\
	&&=\sum_{k=1}^{n}\left(\sum_{i=1}^{k}\sum_{j=1}^{k}a_{kij}+\sum_{i=1}^{k}\sum_{j=1}^{k}a_{ikj}+\sum_{i=1}^{k}\sum_{j=1}^{k}a_{ijk}\right)\nonumber\\
	&&-\sum_{k=1}^{n}\left(\sum_{j=1}^{k}a_{kjk}+\sum_{j=1}^{k}a_{kkj}+\sum_{j=1}^{k}a_{ikk}\right)+\sum_{k=1}^{n}a_{kkk}.\nonumber\\
\end{eqnarray}
In particular, by putting $a_{ijk}=f(i,j)$ in the identity (\ref{ThreeDimKaraji}), we
get

\begin{eqnarray}\label{SpecialtwoDim}
	\sum_{k=0}^{n}\left(\sum_{i=0}^{k}\sum_{j=0}^{k}f(i,j)+\sum_{i=0}^{k}k
	f(k,j)-k f(k,k)\right)=(n+1)\sum_{i=0}^{n}\sum_{i=0}^{n}f(i,j).
\end{eqnarray}
Here $f(i,j)$ is an arbitrary bivaraiate discrete function. 

\begin{example}
	
Put $f(n,k)={n\choose k}$ in formula (\ref{SpecialtwoDim}). Then after some
simplifications, we obtain
$$
\sum_{k=0}^{n}\sum_{j=0}^{k}k{n\choose k}=3n^22^{n-3}+5n2^{n-3}.
$$
	
\end{example}

\begin{example}
	
	Let $a_{ijk}=\frac{1}{i^s}\times\frac{1}{j^s}\times\frac{1}{k^s}$ in
	Formula \ref{ThreeDimKaraji}. Then, we get
	$$\sum_{k=1}^{\infty}\left(\frac{\zeta_{k}(s)}{k^s}-\frac{\zeta_{k}(s)}{k^{2s}}\right)=\frac{\zeta^{2}(s)-\zeta(3s)}{2}.$$

\end{example}

\section{Gosper Decision Algorithm}

We call a sequence of real numbers $\{t_{n}\}_{n\geq 1}$  
\emph{hypergeometric} in variable $n$, if the quotient 
$
\frac{t_{n+1}}{t_{n}}
$
is a \emph{rational} function of $n$. 
More precisely, if there exists two real  
polynomials $p(n)$ and $q(n)$ such that 

\begin{equation}
	\frac{t_{n+1}}{t_{n}} = \frac{p(n)}{q(n)}. 
\end{equation}

\begin{example}
	Let $t_{n} = n!$ for $n\geq 1$. Then, clearly we have 
	$
	\frac{t_{n+1}}{t_{n}} = n+1
	$
	.
	That is, $n!$ is a hypergeometric sequence in $n$. 
	
\end{example}

\begin{example}
	Let $t_{n} = {n \choose k}$, for $1 \leq k \leq n$. 
	Then, we have 
	$$
	\frac{t_{n+1}}{t_{n}} = \frac{{n+1 \choose k}}{ {n \choose k}} = 
	\frac
	{
		\frac
		{(n+1)n!}
		{k! (n+1-k)(n-k)!}	
	}
	{
		\frac{n!}
		{k!(n-k)!}
	}
	=
	\frac{n+1}{n+1-k}
	$$
	Hence, ${n \choose k}$ is a hypergometric sequence in varaiable $n$. 
	One can also simialrly prove that ${n \choose k}$ is a hypergeometric sequence in 
	$k$. 	
\end{example}

\begin{remark}
	
	It is worth to note that the sequence 	
	$
	a_{n,k} = n^{k}
	$
	, 
	for natural numebrs $n$ and $k$, is a hypergeometric sequence in $k$. But, it is not hypergeometric in 
	$
	n
	$
	.
	
\end{remark}
In \cite{Gosp78}, Gosper gives a \emph{decision algorithm} for obtaining the \emph{indefinite sum}
$\sum_{k=1}^{n}f_{k}$. 
More precisely, this algorithm says that  
if $f_{k}$ is a \emph{hypergeometric} sequence, 
then one can find another \emph{hypergeometric} sequence $F_n$ such that
$f_{k}=F_{k+1}-F_{k}$, and therefore after applying the \emph{creative telescoping} \cite{Zeil1991}, 
we find the answer, that is

$$
\sum_{k=1}^{n}f_{k}=\sum_{k=1}^{n}\left\{\left(F_{k}-F_0\right)-\left(F_{k-1}-F_0\right)\right\}=F_{n}. 
$$
The basic idea behind the \emph{Gosper's alogorithm} is the following. 
We already know that 
$
f_{k}=F_{k+1}-F_{k}
$
. 
Now, let us assume that 
$
\frac{F_{k+1}}{F_{k}}
$
is a \emph{rational} function ( since $F_{k}$ is hypergoemtric). 

Hence, we can rewrite the quotient 
$
\frac{f_{k+1}}{f_{k}}
$,
as follows
\begin{eqnarray}
	\frac{f_{k+1}}{f_{k}}
	& = & 
	\frac
	{F_{k+1} - F_{k}}
	{F_{k} - F_{k-1}}, \nonumber\\
	& = & 
	\frac
	{
		\frac{F_{k+1}}{F_{k}} -1
	}
	{
		1 - \frac{F_{k-1}}{F_{k}}
	} 
	.
\end{eqnarray}
That is, the quotient 
$
\frac{f_{k+1}}{f_{k}}
$
is a rational function and therefore $f_{k}$ is also  hypergeometric. 
Now, based on the next lemma, the quotient can always be written in a special form 
which is the \emph{key} for implementing the \emph{Gosper}'s algorithm. 

\begin{lemma}[See \cite{Gosp78}]\label{keyLemmaGos}
	
	A rational function 
	$
	\frac{a(n)}{b(n)}
	$
	can be always written as
	\begin{equation}
		\frac{a(n)}{b(n)} = 
		\frac{p(n) q(n)}{p(n-1) r(n)},
	\end{equation}
	where $p$, $q$ and $r$ are 
	polynomials in $n$ and 
$$
\gcd(q(n), r(n+j) )=1,~~~ \forall j \geq 0.   
$$

\end{lemma} 
Thus, we can prove the following theorem based on Lemma \ref{keyLemmaGos} which guarantees 
the correctness of the Gosper algorithm. 

\begin{theorem}
	
	Let 
	$
	\frac{t(n)}{t(n-1)}
	$
	be a rational function as above. 
	Then
	
	$$
	f(n) = t(n) \frac{p(n)}{q(n+1) t(n)},
	$$
	is a polynomial for which 
	$
	p(n) = q(n+1) f(n) - r(n) f(n-1)
	$
	.

\end{theorem}

\section{Gosper-Karaji Algebraic Identites}

Consider an indefinite sum $\sum_{k=1}^{n}f_{k}$. As we mentioned in the previous section, 
for evaluating this sum, there is a beautiful decision algorithm, Gosper's algorithm
\cite{Gosp78}, which says that if $f_{k}$ is a hypergeometric sequence,
then one can find another hypergeometric sequence $F_n$ such that
$f_{k}=F_{k+1}-F_{k}$, and therefore after applying creative
telescoping, we find the answer, that is
$$
\sum_{k=1}^{n}f_{k}=\sum_{k=1}^{n}\left\{\left(F_{k}-F_0\right)-\left(F_{k-1}-F_0\right)\right\}=F_{n}
$$

Note that we will implicitly assume that 
$
F_{0} = 0
$
.
Otherwise, we will replace the antiderivative 
$
F_{k}
$
,
with the corrected one; that is,
$
F_{k} - F_{0}
$
.
\\
Now we show that for any such indefinite sum, we can obtain a series of identities associated with that identity, using the Karaji's
L-summing method. From now on, we will call them 
\emph{Gosper-Karaji's identities}. 
To do this, we need to apply our general L-summing
equation to an special class of arrays which is called the
\emph{mutiplication table} \cite{H3}. 
Hence, we put
$a_{ij}=t_{i}\times t_{j}$, in formula (\ref{KarajiEq1}). 

After some simple computations, we obtain the following new rearrangement formula

$$
\sum_{k=1}^{n}t_{k}R_{k}=s^{2}_{n},
$$
where $s_{k}=\sum_{i=1}^{k}t_{i}$ and 
$ R_{k}=2s_{k}-t_{k}$. 
Also by the convention, we will assume that $s_{0}=0$.
\\
Now, let's take a closer look at our 
formula. We observe that
$t_{k}R_{k}=(s_{k}-s_{k-1})(s_{k}+s_{k-1})=s^{2}_{k}-s^{2}_{k-1}$.
In fact, we obtain that
$\sum_{k=1}^{n}s^{2}_{k}-s^{2}_{k-1}=s^{2}_{n}$.
By similar
argument, for three dimension (using the formula (\ref{ThreeDimKaraji})), we get 
$$
t_{k}R_{k}=s^{3}_{n}, \hspace{0.5cm}t_{k}R_{k}=3s^{2}_{k}-3s_{k}t_{k}+t^{2}_{k}. 
$$
Indeed, we have $t_{k}R_{k}=s^{3}_{k}-s^{3}_{k-1}$. Therefore, we obtain that
$\sum_{k=1}^{n}s^{3}_{k}-s^{3}_{k-1}=s^{3}_{n}$. 
\\ 
We finally Note that, in general, for
$m$-dimensional multiplication table, we get
$$
t_{k}R_{k}=s^{m}_{n}, \hspace{0.5cm}
t_{k}R_{k}={m\choose 1}s^{m-1}_{k}-{m\choose
	2}s^{m-1}_{k}t_{k}+\cdots+(-1)^{m}{m\choose
	m}t^{m}_{k}
$$
In fact, we get that
$\sum_{k=1}^{n} t_{k}R_{k}  = \sum_{k=1}^{n} (s^{m}_{k}-s^{m}_{k-1})=s^{m}_{n}$. Now we explain our
idea through a very simple but important example.

\begin{example}
	
Consider the simplest indefinite sum $\sum_{k=1}^{n}1$. If we give
	the function $f_{k}=1$ to the latest version of maple package which
	is equipped with Gosper algorithm, we obtain that $s_{k}=k$. Thus,
	by \emph{creative telescoping}, we have
	
	$$
	\sum_{k=1}^{n}1=\sum_{k=1}^{n}\left(k-(k-1)\right)=n.
	$$
	
	Now, Karaji's L-summing method start to do its mission. 
	For $m=2$,
	$m=3$ and $m=4$, we get the following \emph{nontrivial} associated identities, respectively
	$$\sum_{k=1}^{n}\left( 2k-1 \right)=n^2,$$
	$$\sum_{k=1}^{n}{k \choose 2}={n+1 \choose 3},$$
	$$\sum_{k=1}^{n}\left( 2k-1\right)^3={2\sum_{k=1}^{n}\left(2k-1\right) \choose 2}.
	$$

\end{example}

\section{Maple Package Implementation}

The implementation of \emph{Gosper} algorithm in computer algebra system 
\emph{Maple} can be proceed by the software package 
\emph{SumTools}. The package automatically excecute the 
\emph{Gosper's algorithm} and returns the answer if 
$t_{k}$ is a \emph{hypergeomtric} term. 

\begin{example}
	
Let $t_{k} = k$ be a hypergeometric sequence. Then, 
after implementing the library 
\emph{SumTools} of maple:
	
\begin{mapleinput}
	\maplemultiline{> with (sumtools):\\
		[Hypersum; Sumtohyper; extended gosper; gosper; hyperrecursion;\\ hypersum; hyperterm; simpcomb; sumrecursion; sumtohyper]\\
		> sumtools[gosper](k);\\
		\hspace{2cm} \frac{k(k+1)}{2}.  \\}
\end{mapleinput}

Hence, the \emph{Gosper}'s algorithm gives us 
$
s_{k} = \frac{k(k+1)}{2} = {k+1 \choose 2}
$.
Thus, Gosper-Karaji's identites for the cases $m=2,3,4$, 
are as follows

\begin{enumerate}
		\item 
		
		\[
		\sum_{k=1}^{n} {k \choose 1}^{3} = {n+1 \choose 2}^{2}. 
		\]

		\item
		\[
		\frac{3}{4}\sum_{k=1}^{n}{k \choose 1}^{5} 
		+ 
		\frac{1}{4} \sum_{k=1}^{n} {k \choose 1}^{3}
		= 
		{n+1 \choose 2}^{3}
		. 
		\]

		\item
		\[
		\frac{1}{2}\sum_{k=1}^{n}{k \choose 1}^{7} 
		+ 
		\frac{1}{2} \sum_{k=1}^{n} {k \choose 1}^{5}
		= 
		{n+1 \choose 2}^{4}
		. 
		\]

	\end{enumerate}

\end{example}

In the next example, we obtain identites related to the harmonic numbers 
$
H_{n} = \sum_{k=1}^{n} \frac{1}{k}
$
.

\begin{example}
	
	Let $t_{k} = \frac{1}{k}$ be a hypergeometric function in the discrete variable 
	$k$. Then, 	the implementing of the library 
	\emph{SumTools} of maple gives the following result:
	\begin{mapleinput}
		\maplemultiline{> with (sumtools):\\
			[Hypersum; Sumtohyper; extended gosper; gosper; hyperrecursion;\\ hypersum; hyperterm; simpcomb; sumrecursion; sumtohyper]\\
			> sumtools[gosper](\frac{1}{k});\\
			\hspace{2cm} FAIL\\}
	\end{mapleinput}

But, we still can find a solution to the difference equation 
$
t_{k} = s_{k} - s_{k-1}
$ 
based on the definition of \emph{Harmonic} numbers. That is, we get 
$
s_{k} = H_{k}
$
.
Thus, Gosper-Karaji identites for the cases $m=1,2,3$, are as follows
\begin{enumerate}

	\item 
	
	\[
	\sum_{k=1}^{n} \frac{1}{k} = H_{n},
	\]

	\item
	
	\[
	\sum_{k=1}^{n} \frac{1}{k} H_{k} = 
	\frac{1}{2} H^{2}_{n} + \frac{1}{2} H^{(2)}_{n},
	\]

	\item

	\[
	\sum_{k=1}^{n} \frac{1}{k} H^{2}_{k}  - 
	\sum_{k=1}^{n} \frac{1}{k^2} H_{k}
	= 
	\frac{1}{3} H^{3}_{n} - \frac{1}{3} H^{(3)}_{n}
	.
	\]

\end{enumerate}

\end{example}

\begin{example}
	
Let $t_{k} = \frac{1}{k(k+1)}$ be a hypergeometric function in the discrete variable $k$. 
	Then, the implementing of the library 
	\emph{SumTools} of maple gives the following result:
	
	\begin{mapleinput}
		\maplemultiline{> with (sumtools):\\
			[Hypersum; Sumtohyper; extended gosper; gosper; hyperrecursion;\\ hypersum; hyperterm; simpcomb; sumrecursion; sumtohyper]\\
			> sumtools[gosper](\frac{1}{k(k+1)});\\
			\hspace{2cm} FAIL  \\}
	\end{mapleinput}	
	
Again, the Gosper algorithm fails but still we can use the idea of \emph{creative telescoping} \cite{Zeil1991} 
to find the solution  
$
s_{k} = 2-\frac{k!}{(2k+1)!}
$. 
Thus, Gosper-Karaji identites for the cases $m=1,2,3$, are as follows
	
	\begin{enumerate}

		\item 
		
		\[
		\frac{1}{2}	
		\sum_{k=1}^{n} \frac{1}{{k+1 \choose 2} } = 1-\frac{1}{n+1},
		\]

		\item
		
		\[
		\sum_{k=1}^{n} \frac{1}{{k+1 \choose 1}^{2} } 
		+ 
		\frac{1}{4}	
		\sum_{k=1}^{n} \frac{1}{{k+1 \choose 2}^{2} }
		= 1-\frac{1}{(n+1)^{2}},
		\]

		\item

		\[
		\frac{12}{8} 
		\sum_{k=1}^{n} \frac{
			{k+1\choose2} {k \choose 2} 
		}
		{
			{k+1 \choose 2}^{3}
		} 
		+ 
		\frac{1}{8}
		\sum_{k=1}^{n} \frac{1}
		{
			{k+1 \choose 2}^{3}
		}
		=
		\frac
		{
			{n \choose 1}^{3}
		}
		{
			{n+1 \choose 1}^{3}
		}
		.
		\]

	\end{enumerate}

\end{example}

\begin{example}
	
	Assume  
	$
	t_{k} = (4k+1)\frac{k!}{(2k+1)!}
	$ 
	be a hypergeometric function in the discrete 
	variable $k$. 
	Then, the implementing of the library 
	\emph{SumTools} of maple gives the following result:
	
\begin{mapleinput}
		\maplemultiline{> with (sumtools):\\
			[Hypersum; Sumtohyper; extended gosper; gosper; hyperrecursion;\\ hypersum; hyperterm; simpcomb; sumrecursion; sumtohyper]\\
			> sumtools[gosper](\frac{1}{k(k+1)});\\
			\hspace{2cm} 1-\frac{1}{k+1}  \\}
\end{mapleinput}

Thus, Gosper-Karaji's identites for 
the cases $m=1,2,3$, are as follows
	
	\begin{enumerate}

		\item 
		
		\[
		\sum_{k=0}^{n} 
		(4k+1)\frac{k!}{(2k+1)!}
		= 
		2-\frac{n!}{(2n+1)!},
		\]

		\item

		\[
		\sum_{k=0}^{n} 
		(4k+3)(4k+1)\frac{k!^{2}}{(2k+1)!^{2}}
		= 
		4-\frac{n!^{2}}{(2n+1)!^{2}},
		\]


	\end{enumerate}

\end{example}



\end{document}